\documentclass{amsart}
\usepackage{amsmath, amssymb, latexsym}
\title{A Kruskal-Katona Type Theorem for Graphs}
\author{Andy Frohmader}
\address{Department of Mathematics, University of
  Washington, Seattle, WA 98195-4350}
\email{frohmade@math.washington.edu}

\newtheorem{theorem}{Theorem}[section]
\newtheorem{proposition}[theorem]{Proposition}

\newtheorem{lemma}[theorem]{Lemma}
\newtheorem{definition}[theorem]{Definition}

\newtheorem{example}[theorem]{Example}
\newtheorem{construction}[theorem]{Construction}

\newcommand{\link}[2]{\ensuremath{\textup{lk}_{#1}({#2})}}
\newcommand{\oldbd}[2]{\ensuremath{\textup{old}_{#1}({#2})}}
\newcommand{\smbd}[2]{\ensuremath{\textup{sm}_{#1}({#2})}}
\newcommand{\lgbd}[2]{\ensuremath{\textup{lg}_{#1}({#2})}}
\newcommand{\conbd}[2]{\ensuremath{\textup{con}_{#1}({#2})}}
\newcommand{\ratio}[2]{\ensuremath{\textup{rat}_{#1}({#2})}}
\newcommand{\clique}[2]{\ensuremath{\textup{cl}_{#1}({#2})}}
\newcommand{\spand}{\ensuremath{\qquad \mbox{ \textup{and}}}}

\def\proof{\smallskip\noindent {\it Proof: \ }}
\def\endproof{\hfill\ensuremath{\square}\medskip}

\begin{document}

\maketitle

\begin{abstract}
A bound on consecutive clique numbers of graphs is established.  This bound is evaluated and shown to often be much better than the bound of the Kruskal-Katona theorem.  A bound on non-consecutive clique numbers is also proven.
\end{abstract}

\section {Introduction}

The goal of this paper is to provide a partial answer to the question of how many $(k+1)$-cliques a graph with a given number of $k$-cliques can have.

Given a graph, we can readily count its cliques of various sizes.  For simplicity, we refer to a clique on $n$ vertices as an $n$-clique.  We can count the number of cliques of all possible sizes to get the clique vector of a graph. It then makes sense to ask what integer vectors can arise as clique vectors of graphs.

Simplicial complexes can be thought of as generalizations of graphs, as from any graph, we can form its clique complex, a simplicial complex whose faces correspond to the cliques of the graph.  The question of which integer vectors can be clique vectors of graphs then becomes a question of which integer vectors can be face vectors of simplicial complexes.

This question was answered in the 1960s independently by Kruskal \cite{kruskal} and Katona \cite{katona}.  Much additional work on face vectors has been done since then.  For example, Stanley \cite{cohenface} characterized the face vectors of Cohen-Macaulay complexes, while Frankl, F\"{u}redi, and Kalai \cite{balanced} characterized the face vectors of complexes with a given chromatic number.

Simplicial complexes that arise as clique complexes of graphs are also known as flag complexes, and are of interest in their own right in this context.  For example, the face ideal (see \cite{greenbook}, chapter 2) of a flag complex is generated by quadratic monomials.  Thus, the face ring of a simplicial complex is Koszul exactly if the complex is a flag complex (see \cite{koszul}).

Work toward characterizing the clique vectors of graphs goes as far back as Mantel's theorem (see \cite{mantel}, p. 31), which states that a graph with $n$ vertices and no triangles has at most $\big\lfloor{n^2 \over 4}\big\rfloor$ edges.  Tur\'{a}n's theorem \cite{turan} generalized this to give a bound on the most edges a graph could have in terms of the number of vertices and size of the largest clique. Zykov \cite{zykov} further generalized this to get a bound on cliques of all sizes in terms of the number of vertices and size of the largest clique.  More recently, Eckhoff \cite{eckhoff0, eckhoff} proved bounds on clique numbers in terms of the number of edges and size of the largest clique.  A paper of the author \cite{previous} generalized these results to give a bound on the number of $(i+1)$-cliques of a graph in terms of the number of $i$-cliques and the size of the largest clique of the graph, thereby affirming a conjecture of Kalai (unpublished; see \cite[p.~100]{greenbook}) and Eckhoff \cite{mainconj}.

A related result of Lov\'{a}sz and Simonovits \cite{lovasz} proved a lower bound on the number of $k$-cliques a graph with a given number of vertices and edges must contain.  A very recent paper of Herzog, Hibi, Murai, Trung, and Zhang \cite{hhmtz} characterized the clique vectors of chordal and strongly chordal graphs.

However, the problem of characterizing the clique vectors of general graphs remains open.  Clique vectors of graphs must satisfy the inequalities shown by Kruskal and Katona to characterize face vectors of all simplicial complexes.  Satisfying these inequalities is enough to ensure that an integer vector is the face vector of a simplicial complex, but not necessarily the clique vector of a graph.  For example, there is a complex with 9 faces on 3 vertices and 3 faces on 4 vertices, but it is easy to show by checking cases that there is no graph with exactly 9 3-cliques and 3 4-cliques.

Dealing with graphs is harder than simplicial complexes because we can readily add a single face of arbitrary dimension to a simplicial complex, but can only directly control the vertices and edges of graphs.  Adding a single edge often adds many cliques of each of several sizes to a graph, meaning that a single construction is not enough as in the case of simplicial complexes.

For example, a complete graph on 10 vertices contains 120 3-cliques.  If one edge is removed, the resulting graph contains 112 3-cliques. Both of these graphs attain the bound of the Kruskal-Katona theorem for the most 4-cliques a graph can have in terms of the number of 3-cliques, but they shed no light on how many 4-cliques a graph with 119 3-cliques can have, and there is not a natural intermediate construction that does.

As another example, a complete graph on 7 vertices has 35 3-cliques and 35 4-cliques.  However, by \cite{previous}, if a graph has 35 3-cliques and does not have a 5-clique, then it can have at most 17 4-cliques.  This result can give a useful bound on graphs that do not have a large clique, but if we allow the graph to have large enough cliques, the result of \cite{previous} usually coincides with the bounds of the Kruskal-Katona theorem.

The structure of the paper is as follows.  In Section 2, we give precise definitions of some needed concepts, including some mentioned above.  We also state some theorems needed as background material.  The section concludes by stating our main result, Theorem~\ref{maintheorem}.

The main focus of Section 3 is to address what happens to graphs that do have a large clique.  We derive a bound in Lemma~\ref{lemma5} on the clique numbers of graphs that do have as large of a clique as is possible without exceeding the allowed number of cliques.  For example, if we are given that a graph has 102 3-cliques, then the graph could contain a clique on 9 vertices, as a 9-clique only contains ${9 \choose 3} = 84$ 3-cliques. However, it could not contain a clique on 10 vertices, as such a clique contains ${10 \choose 3} = 120$ 3-cliques, which is more than allowed.  The main result of this section states that a graph with 102 3-cliques and at least one 9-clique can contain at most 147 4-cliques.

Continuing this example, \cite{previous} states that a graph with 102 3-cliques and no 9-cliques can have at most 146 4-cliques.  Since any graph with 102 3-cliques either contains a 9-clique or it does not, its number of 4-cliques must be bounded above by either 147 or 146.  Hence, it is bounded above by the larger value, 147.  For comparison, the Kruskal-Katona theorem states that the graph can have at most 149 4-cliques.

In Section 4, we evaluate our bounds to see just how good they are compared to the Kruskal-Katona theorem.  Theorem~\ref{convprob} is a convergence in probability result which shows that when a graph does have the largest clique possible, the number of $(k+1)$-cliques allowed by the bound of Theorem~\ref{maintheorem} is nearly always much closer to the most $(k+1)$-cliques of any graph with the chosen number of $k$-cliques than to the bound of the Kruskal-Katona theorem.  Proposition~\ref{smallstrict} shows that the bound when we do not have a large clique is always strictly smaller than the bound of the Kruskal-Katona theorem.

We give three constructions in Section 4.3 to show that our bounds are attained by graphs under certain circumstances.  These constructions are the motivation behind the particular bounds that are proven.  To finish the running example, Construction~\ref{const2} provides a graph with 102 3-cliques and 147 4-cliques, so the bound of Theorem~\ref{maintheorem} is attained.

Finally, in Section 5, we consider what happens to non-consecutive clique numbers of graphs, such as how many 7-cliques a graph can have in terms of its number of 4-cliques.  A bound in Theorem~\ref{nonconsec} immediately follows from work earlier in the paper, just as bounds on non-consecutive face numbers often followed immediately from bounds on consecutive face numbers in previous work on face numbers by various authors.

Unlike when dealing with simplicial complexes, however, our bound on non-consecutive face numbers surprisingly gives new information. We demonstrate this in Example~\ref{example} to show that bounds on consecutive clique numbers are not enough to characterize the clique vectors of graphs.

Throughout this paper, most of the lemmas are elementary algebraic statements, but the proofs are often combinatorial, and involve constructing various graphs or simplicial complexes.  That the Kruskal-Katona theorem or various other related results apply to the construction is used extensively in proving the needed results.

\section {Background and definitions}

In this section, we review some material that will be needed for our results.

\subsection {Graphs and simplicial complexes}

Recall that a \textit{graph} $G$ is a set $V$ of vertices and a set $E$ of edges connecting pairs of vertices.  This paper deals only with simple graphs on a finite vertex set without loops or multiple edges. A \textit{clique} of a graph is a complete subgraph, that is, a subset $C \subset V$ of the vertices such that every two vertices of $C$ are connected by an edge.  In particular, if $C$ contains only one vertex, it is a clique without any condition on edges.  Every graph also has a unique clique on zero vertices. We can count the number of cliques of a given size.

\begin{definition}
\textup{The \textit{$i$-th clique number} of a graph $G$, denoted \clique{i}{G}, is the number of cliques of $i$ vertices in $G$. These are also called \textit{$i$-cliques} of $G$. If the largest clique of $G$ has $d$ vertices, the \textit{clique vector} of $G$ is the vector
$$\clique{}{G} = (\clique{0}{G}, \clique{1}{G}, \dots , \clique{d}{G}).$$}
\end{definition}

While the main theorems are results about cliques of graphs, the proofs extensively use simplicial complexes.  Recall that a \textit{simplicial complex} $\Delta$ on a vertex set $V$ is a collection of subsets of $V$ such that (i) for every $v \in V$, $\{v\} \in \Delta$ and (ii) for every $B \in \Delta$, if $A \subset B$, then $A \in \Delta$.  The elements of $\Delta$ are called \textit{faces}.  A face on $i$ vertices is said to have \textit{dimension} $i-1$, while the dimension of a complex is maximum dimension of a face of the complex.  The maximal faces (under inclusion) are called \textit{facets}.  A simplicial complex in which all maximal faces are of the same dimension is called \textit{pure}.

We can count the number of faces on a given number of vertices in a simplicial complex, just as we can count cliques in graphs.

\begin{definition}
\textup{The \textit{$i$-th face number} of a simplicial complex $C$, denoted \clique{i}{C} is the number of faces in $C$ containing $i$ vertices. These are also called \textit{$i$-faces} of $C$. If dim $C = d-1$, the \textit{face vector} of $C$ is the vector
$$\clique{}{C} = (\clique{0}{C}, \clique{1}{C}, \dots , \clique{d}{C}).$$}
\end{definition}

It is sometimes useful in inductive proofs to consider certain subcomplexes of a given simplicial complex, such as its links.

\begin{definition}
\textup{Let $\Delta$ be a simplicial complex and $F \in \Delta$.  The \textit{link} of $F$, \link{\Delta}{F}, is defined as}
$$\link{\Delta}{F} := \{G \in \Delta \ |\ F \cap G = \emptyset, F \cup G \in \Delta\}.$$
\end{definition}

The link of a face of a simplicial complex is itself a simplicial complex.  We can analogously define a link of a clique of a graph.

\begin{definition}
\textup{The link of a clique $C = \{v_1, v_2, \dots, v_n\}$ of vertices in a graph $G$, denoted \link{G}{v_1v_2\dots v_n}, is the induced subgraph of $G$ on the set of vertices that are adjacent to all vertices of $C$.}
\end{definition}

This paper usually considers only the link of a single vertex, or at most, the link of an edge.

A useful construction in building certain simplicial complexes is the reverse-lexicographic (``rev-lex") order. To define the rev-lex order of $i$-faces of a simplicial complex on $n$ vertices, we start by labelling the vertices $1, 2, \dots$. Let $\mathbb{N}$ be the set of natural numbers, let $A$ and $B$ be distinct subsets of $\mathbb{N}$ with $|A| = |B| = i$, and let $A \nabla B$ be the symmetric difference of $A$ and $B$.

\begin{definition}
\textup{For $A, B \subset \mathbb{N}$ with $|A| = |B|$, we say that $A$ precedes $B$ in the rev-lex order if max$(A \nabla B) \in B$, and $B$ precedes $A$ otherwise.}
\end{definition}

For example, $\{2, 3, 5\}$ precedes $\{1, 4, 5\}$, as 3 is less than 4, and $\{3, 4, 5\}$ precedes $\{1, 2, 6\}$.

\begin{definition}
\textup{The \textit{rev-lex complex on $m$ $i$-faces} is the pure complex whose facets are the first $m$ $i$-sets possible in rev-lex order.  This complex is denoted $C_i(m)$.}
\end{definition}

We can also specify more than one number in the face vector.  For two sequences $i_1 < \dots < i_r$ and $(m_1, \dots, m_r)$, let
$$C = C_{i_1}(m_1) \cup C_{i_2}(m_2) \cup \dots \cup C_{i_r}(m_r).$$
A standard way to prove the Kruskal-Katona theorem involves showing that if the numbers $m_1, \dots, m_r$ satisfy the bounds of the theorem, then the complex $C$ has exactly $m_j$ $i_j$-faces for all $j \leq r$ and no more.  In this case, we refer to $C$ as the \textit{rev-lex complex on $m_1$ $i_1$-faces, \dots, $m_r$ $i_r$-faces}.

The notion of rev-lex complexes can be extended to colored complexes.  The \textit{chromatic number} of a simplicial complex is the minimal number of colors required to color all vertices of the complex such that no two vertices in any face are the same color.  This definition coincides with the chromatic number of the 1-skeleton of the complex, taken as a graph.

\begin{definition}
\textup{A subset $A \subset \mathbb{N}$ is \textit{$r$-permissible} if, for every two $a, b \in A$, $r$ does not divide $a-b$.  The \textit{$r$-colored rev-lex complex on $m$ $i$-faces} is the pure complex whose facets are the first $m$ $r$-permissible $i$-sets in rev-lex order.}
\end{definition}

We can specify more than one number in the face vector for colored rev-lex complexes in the same manner as for the usual (uncolored) rev-lex complexes.

\subsection {Lemmas on binomial representations}

In this section, we give some basic lemmas, which are necessary in order for the bounds on clique numbers to be well-defined.  We start with some notation.

\begin{definition}
For integers $k \geq s \geq 0$, define
$$r_k(n_k, n_{n-1}, \dots , n_{k-s}) = {n_k \choose k} + {n_{k-1} \choose k-1} + \dots + {n_{k-s} \choose k-s}.$$
\end{definition}

For example, the basic identity ${n \choose k} = {n-1 \choose k} + {n-1 \choose k-1}$ can be expressed as $r_k(n) = r_k(n-1, n-1)$. This expression is not unique, as the same identity can be expressed as $r_k(n) = r_k(n-1) + r_{k-1}(n-1)$, $r_{k+1}(0, n) = r_{k+1}(0, n-1) + r_{k-1}(n-1)$, or in many other ways.  We can, however, make it unique with additional restrictions.

\begin{lemma}\label{kklemma}
Given positive integers $m$ and $k$, there are unique integers $s \geq 0$ and $n_k > n_{k-1} > \dots > n_{k-s} \geq k-s > 0$ such that $m = r_k(n_k, n_{k-1}, \dots, n_{k-s})$.
\end{lemma}

This is a standard lemma associated with the Kruskal-Katona theorem \cite{katona, kruskal}, so we do not give a proof here.

One convention we use throughout this paper is that any time we define constants $a_k, a_{k-1}, \dots , a_{k-s}$ by saying that $r_k(a_k, a_{k-1}, \dots, a_{k-s})$ is equal to a particular constant, the $a_i$s are the unique choice of constants that satisfy the conditions of Lemma~\ref{kklemma}.  In particular, if $r_k(a_k, a_{k-1}, \dots, a_{k-s})$ appears in the statement of a lemma and is the first time that the $a_i$s have appeared, they are defined to be the unique constants satisfying Lemma~\ref{kklemma} to make $r_k(a_k, a_{k-1}, \dots, a_{k-s})$ equal to some particular constant.  This convention only applies when we are defining new constants, and not merely using constants that were previously defined in the proof.

The value of $s$ often does not matter to the proof.  For notational simplicity, we often leave off the last term and talk of $a_k, a_{k-1}, \dots$.

\begin{lemma}\label{lgbdlemma}
Given positive integers $m$ and $k$, there are unique integers $s \geq 0$, $n_k > n_{k-1} \geq k-2$, and $a_{k-1} > a_{k-2} > \dots > a_{k-s} \geq k-s > 0$, such that
\begin{eqnarray*}
r_{k-2}(n_{k-1}) & > & r_{k-1}(a_{k-1}, \dots, a_{k-s}) \spand \\ m & = & r_k(n_k, n_{k-1}) + r_{k-1}(a_{k-1}, \dots, a_{k-s}).
\end{eqnarray*}
\end{lemma}

A bit of interpretation is required here for the case $s = 0$.  This corresponds to the case when $m = r_k(n_k, n_{k-1})$, and the other conditions on the $a_i$s are considered to be trivially satisfied. Similarly, the case $n_{k-1} = k-2$ corresponds to the case when $m = r_k(n_k)$.

\proof  Define $n_k$ and $n_{k-1}$ such that $m = r_k(n_k, n_{k-1}, n_{k-2}, \dots)$ is the unique representation of Lemma~\ref{kklemma}.  Let $q = m - r_k(n_k, n_{k-1})$, and define $a_i$s such that $q = r_{k-1}(a_{k-1}, a_{k-2}, \dots, a_{k-s})$.  We must have $r_{k-2}(n_{k-1}) > q$ or else $n_{k-1}$ would have been chosen to be larger, so this satisfies the conditions of the lemma.

For uniqueness, once we pick $n_k$ and $n_{k-1}$, the $a_i$s are forced to be unique.  If we make $n_k$ one larger, then $r_k(n_k) > m$.  If we make $n_{k-1}$ one larger, then $r_k(n_k, n_{k-1}) > m$.  If we make $n_k$ or $n_{k-1}$ smaller, then we reduce $r_k(n_k, n_{k-1})$ by at least $r_{k-2}(n_{k-1})$, which would force $q \geq r_{k-2}(n_{k-1})$.  In any of these cases, it is not possible to pick $a_i$s to satisfy the lemma, so the choices of $n_k$ and $n_{k-1}$ are also unique.  \endproof

While there are numbers put into the format of this lemma at various places, it is explicitly stated when conditions beyond those of Lemma~\ref{kklemma} are assumed to be satisfied.

\begin{definition}
\textup{The Tur\'{a}n graph $T_{n,r}$ is the graph obtained by partitioning $n$ vertices into $r$ parts as evenly as possible, and making two vertices adjacent exactly if they are not in the same part.  Define ${n \choose k}_r := \clique{k}{T_{n,r}}$.}
\end{definition}

\begin{lemma}\label{colorcan}
Given positive integers m, k, and r with $r\geq k$, there are unique $s$, $n_k$, $n_{k-1}$, \dots , $n_{k-s}$ such that
$$m = {n_k\choose k}_r + {n_{k-1}\choose k-1}_{r-1} + \dots + {n_{k-s}\choose k-s}_{r-s},$$
$n_{k-i}-\big\lfloor{n_{k-i}\over r-i}\big\rfloor > n_{k-i-1}$ for all $0\leq i < s,$ and $n_{k-s}\geq k-s > 0$.
\end{lemma}

The original use of this lemma in \cite{balanced} misstated it. A correct version that is equivalent to the above lemma appears in \cite[Theorem 15.1.3]{colorfixed}.

\subsection {Kruskal-Katona type theorems}

We need some notation to simplify the discussion of the bounds to be proven.

\begin{definition}
\textup{Let $m = r_k(n_k, n_{k-1})+r_{k-1}(a_{k-1}, a_{k-2}, \dots)$ be the representation of Lemma~\ref{lgbdlemma}. Define
$$\lgbd{k}{m} := r_{k+1}(n_k, n_{k-1})+r_k(a_{k-1}, a_{k-2}, \dots).$$}
\end{definition}

\begin{definition}
\textup{Let $m = r_k(n_k, n_{k-1}, \dots)$.  Define $$\oldbd{k}{m} := r_{k+1}(n_k, n_{k-1}, \dots).$$
If $n_k > k$, then let $a_k, a_{k-1}, \dots$ and $s$ be the unique integers satisfying the conditions of Lemma~\ref{colorcan} such that
$$m = {a_k\choose k}_{n_k-1} + {a_{k-1}\choose k-1}_{n_k-2} + \dots + {a_{k-s}\choose k-s}_{n_k-s-1}.$$
Define
$$\smbd{k}{m} := {a_k\choose k+1}_{n_k-1} + {a_{k-1}\choose k}_{n_k-2} + \dots + {a_{k-s}\choose k-s+1}_{n_k-s-1}.$$}
\end{definition}

If $n_k = k$, then $\smbd{k}{m}$ is undefined; in this case, $\oldbd{k}{m} = 0$.

The \oldbd{k}{m} and \smbd{k}{m} bounds have already been proven in the relevant cases, so we merely cite them here.

\begin{theorem}
[Kruskal-Katona \cite{katona, kruskal}] Let $C$ be a simplicial complex. If $\clique{k}{C} = m$, then $\clique{k+1}{C} \leq \oldbd{k}{m}$. Furthermore, if a non-negative integer vector $f =(1, c_1, c_2, \dots)$ satisfies these inequalities for all $k$, then there is a rev-lex complex $C$ with $f$ as its face vector.
\end{theorem}

\begin{theorem}
[Frankl-F\"{u}redi-Kalai \cite{balanced}] \label{coloredkk} For an $r$-colorable simplicial complex $C$, let
$$m = \clique{k}{C} = {n_k\choose k}_r + {n_{k-1}\choose k-1}_{r-1} + \dots + {n_{k-s}\choose k-s}_{r-s}$$
be the unique representation of Lemma~\ref{colorcan}. Then
$$\clique{k+1}{C} \leq {n_k\choose k+1}_r + {n_{k-1}\choose k}_{r-1} + \dots + {n_{k-s}\choose k-s+1}_{r-s}.$$
Furthermore, given a vector $f = (1, c_1, c_2, \dots c_t)$ that satisfies this bound for all $1 \leq k < t$, there is an $r$-colorable rev-lex complex that has $f$ as its face vector.
\end{theorem}

If $r_k(r+1) \leq m < r_k(r+2)$, this theorem states that $\clique{k+1}{C} \leq \smbd{k}{m}$.

\begin{theorem}\label{mytheorem}
For a positive integer $r$ and a graph $G$ with $\clique{r+1}{G} = 0$, let
$$m = \clique{k}{G} = {n_k\choose k}_r + {n_{k-1}\choose k-1}_{r-1} + \dots + {n_{k-s}\choose k-s}_{r-s}$$
be the unique representation of Lemma~\ref{colorcan}. Then
$$\clique{k+1}{G} \leq {n_k\choose k+1}_r + {n_{k-1}\choose k}_{r-1} + \dots + {n_{k-s}\choose k-s+1}_{r-s}.$$
\end{theorem}

This theorem was proven by the author in \cite{previous}.  It verified a statement conjectured independently by Kalai (unpublished; see \cite[p.~100]{greenbook}) and Eckhoff \cite{mainconj}.  Together with Theorem~\ref{coloredkk}, it implies that for every clique complex of dimension $r-1$, there is an $r$-colorable complex with the same face vector.

If we let $\clique{k}{G} = r_k(n_k, n_{k-1}, \dots)$ be the representation of Lemma~\ref{kklemma}, then if $G$ does not have an $n_k$-clique, this immediately implies $\clique{k+1}{G} \leq \smbd{k}{m}$.  Furthermore, if $n_k = n_{k-1}+1$, then because ${n_k+1 \choose k}_{n_k} = r_k(n_k, n_k-1)$, this theorem states that $\clique{k}{G} \leq \lgbd{k}{G}$.  Lemma~\ref{lemma5} shows that  $\clique{n_k}{G} > 0$ is also a sufficient condition for $\clique{k}{G} \leq \lgbd{k}{G}$.  Combining these results gives the statement of our main theorem.

\begin{theorem}\label{maintheorem}
Let $G$ be a graph and $\clique{k}{G} = m$.  Then
$$\clique{k+1}{G} \leq \max\{\lgbd{k}{m}, \smbd{k}{m}\}.$$
\end{theorem}

The notation \smbd{k}{m} is chosen because it is the bound that applies when the largest clique of the graph is ``small", that is, not as large as it could have been for the allowed number of $k$-cliques.  Likewise, \lgbd{k}{m} was so named because it is the bound that applies when the largest clique is as large as it possibly could have been.  The notation \oldbd{k}{m} is used because that bound is much older than the others, having been first proven in the 1960s.

The goal of the next section is to prove Theorem~\ref{maintheorem}.

\section {Proof of the main theorem}

Theorem~\ref{mytheorem} gives a bound on clique numbers of graphs that depends on the size of the largest clique of the graph.  If the largest clique of the graph is relatively small, this bound can be much less than the bound of the Kruskal-Katona theorem, and allows far fewer $(k+1)$-cliques than a graph with a larger clique can be readily constructed to have.  However, if a graph with a prescribed number of $k$-cliques has the largest clique it could possibly have without exceeding the allowed number of $k$-cliques, this bound often coincides with the Kruskal-Katona theorem.  In this section, we prove a bound on the number of $(k+1)$-cliques that such a graph can have.

We need several technical lemmas.  The lemmas are stated in terms of elementary algebra, though their proofs are often combinatorial and involve constructing simplicial complexes.  The lemmas lead to Lemma~\ref{lemma5}, which is a bound on clique numbers that applies to graphs that do have the largest clique possible.  Our main theorem then follows from a combination of Lemma~\ref{lemma5} and Theorem~\ref{mytheorem}.

\begin{lemma}\label{dimshift}
If $j > 0$, $k > 0$, and $r_k(a_k, a_{k-1}, \dots) \geq r_k(b_k, b_{k-1}, \dots),$ then $r_j(a_k, \dots) \geq r_j(b_k, \dots)$.
\end{lemma}

\proof  If $a_i$ exists and $b_i$ does not, then we will use the convention that $a_i > b_i$.  If $a_i = b_i$ for all $i$, the result is trivial.  Otherwise, let $m = \max\{i\ |\ a_i \not = b_i\}$.  If $a_m$ does not exist, then $r_k(a_k, a_{k-1}, \dots) < r_k(b_k, b_{k-1}, \dots),$ a contradiction. Otherwise, we can subtract $r_k(a_k, \dots, a_{m+1})$ from both sides of the statement of the lemma to get $r_m(a_m, \dots) \geq r_m(b_m, \dots)$.  If $b_m > a_m$, then $b_m \geq a_m+1$, so
$$r_m(b_m, \dots) \geq r_m(b_m) \geq r_m(a_m+1) > r_m(a_m, \dots),$$
a contradiction.  Thus, $a_m > b_m$, and so
$$\hspace{43 pt} r_j(a_k, \dots) \geq r_j(a_k, \dots, a_m) \geq r_j(b_k, \dots, b_m+1) \geq r_j(b_k, \dots). \hspace{43 pt} \square$$

Lemma~\ref{dimshift} can also be derived from the Kruskal-Katona theorem with a comparably easy proof.

\begin{lemma} \label{disjoint}
If $m = r_k(c_k, c_{k-1}, \dots) = r_k(a_k, a_{k-1}, \dots) + r_k(b_k, b_{k-1}, \dots)$, then
$$r_{k+1}(c_k, \dots) \geq r_{k+1}(a_k, \dots) + r_{k+1}(b_k, \dots).$$
\end{lemma}

\proof  The rev-lex complex on $r_{k+1}(a_k, \dots)$ $(k+1)$-faces has $r_k(a_k, \dots)$ $k$-faces.  The rev-lex complex on $r_{k+1}(b_k, \dots)$ $(k+1)$-faces has $r_k(b_k, \dots)$ $k$-faces.  Then their disjoint union has $m$ $k$-faces and $r_{k+1}(a_k, \dots) + r_{k+1}(b_k, \dots)$ $(k+1)$-faces. By the Kruskal-Katona theorem, if a simplicial complex has $m$ $k$-faces, then it has at most $\oldbd{k}{m} = r_{k+1}(c_k, \dots)$ $(k+1)$-faces, so the complex constructed above as the disjoint union of two others satisfies this bound.  \endproof

An equivalent formulation of the above lemma is $\oldbd{k}{m+n} \geq \oldbd{k}{m} + \oldbd{k}{n}$.

\begin{lemma} \label{linksub}
If $r_k(a_k, a_{k-1}, \dots) \geq r_k(b_k, b_{k-1}, \dots)$ and $r_k(c_k, c_{k-1}, \dots) =$ \newline $r_k(a_k, a_{k-1}, \dots) + r_{k-1}(b_k, b_{k-1}, \dots)$, then
$$r_{k+1}(c_k, \dots) \geq r_{k+1}(a_k, \dots) + r_k(b_k, \dots).$$
\end{lemma}

\proof  Let $C$ be the rev-lex complex on $r_{k+1}(a_k, \dots)$ $(k+1)$-faces, and let $D$ be the rev-lex complex on $r_{k+1}(b_k, \dots)$ $(k+1)$-faces.  Since $r_k(a_k, \dots) \geq r_k(b_k, \dots)$, by Lemma~\ref{dimshift}, $C \supseteq D$. Form a new complex $E$ by taking $C$ and adding a new vertex $v$, such that $\link{E}{v} = D$.

The number of $k$-faces of $E$ is the number of $k$-faces containing $v$, plus the number not containing $v$.  These are \clique{k}{C} and \clique{k-1}{D}, respectively, so
$$\clique{k}{E} = \clique{k}{C} + \clique{k-1}{D} = r_k(a_k, \dots) + r_{k-1}(b_k, \dots) = r_k(c_k, \dots).$$
By the Kruskal-Katona theorem, $\clique{k+1}{E} \leq r_{k+1}(c_k, \dots)$.  Applying the same argument for the number of $(k+1)$-faces of $E$ gives
$$\hspace{18 pt} r_{k+1}(c_k, \dots) \geq \clique{k+1}{E} = \clique{k+1}{C} + \clique{k}{D} = r_{k+1}(a_k, \dots) + r_k(b_k, \dots). \hspace{18 pt} \square$$

The next lemma has an algorithmic proof, and is used repeatedly in this paper, both in the proof of the main theorem and later.

\begin{lemma} \label{algorithm}
If $m = r_k(c_k, c_{k-1}, \dots) + r_k(a_k, a_{k-1}, \dots)$, $c_k \leq a_k$, and $m = r_k(a_k+1) + r_k(b_k, b_{k-1}, \dots)$, then
$$r_{k+1}(a_k+1) + r_{k+1}(b_k, \dots) > r_{k+1}(c_k, \dots) + r_{k+1}(a_k, \dots).$$
\end{lemma}

\proof  Set up a rectangular board with two rows and $k$ columns.  In each square of the board, we can either write a positive integer or leave the square blank.  Number the columns based on how far from the right edge they are.  The far right column is column 1, then one next to it is column 2, and so forth, with column $k$ being the far left one.  An arrangement of numbers on the board is \textit{permissible} if \begin{enumerate}
\item for each pair of adjacent squares in the same row, either the one to the left contains a larger number than the one to the right or the one on the right is empty and
\item for each row, if the rightmost non-empty column in a row is column $i$, then the entry in that box is at least $i$.
\end{enumerate}

Let the numbers in the top row be $x_k, x_{k-1}, \dots, x_g$ and the numbers in the bottom row be $y_k, y_{k-1}, \dots, y_h$. A rearrangement of the numbers on the board (or a \textit{move}) is \textit{allowable} if
\begin{enumerate}
\item the arrangement of numbers on the board after the move is permissible,
\item the sum $r_k(x_k, \dots, x_g) + r_k(y_k, \dots, y_h)$ is unchanged,
\item the sum $r_{k+1}(x_k, \dots, x_g) + r_{k+1}(y_k, \dots, y_h)$ does not decrease, and
\item $r_k(x_k, \dots, x_g)$ strictly increases.
\end{enumerate}

The structure of the proof is to have the board start with $a_k, a_{k-1}, \dots$ as the entries in the top row and $c_k, c_{k-1}, \dots$ as the entries in the bottom row, with any leftover boxes initially empty.  We then define a number of moves that are allowable under certain circumstances and show that in all possible circumstances, there is an allowable move, until the board reaches a state in which the top row has $a_k+1$ in column $k$ and the rest of the row is empty, while the bottom row has $b_k, b_{k-1}, \dots$ as its entries.  Conditions two, three, and four are usually trivial to check, so we do not give reasons why they hold in such cases.  The result of the lemma follows from the conditions for a move to be allowable and that at least one move strictly increases the sum of condition 3.

We now explain the needed types of allowable moves.  Figure~1 contains a flow chart showing how to choose which move to make at a given step.

Suppose that $g > h$ and $x_g \geq y_g$.  A move of the first type is to move the last $g-h$ entries in the bottom row up to the top row.  This gives a permissible arrangement of the board since $x_g \geq y_g > y_{g-1}$.

Suppose that there is an $i$ for which $y_i > x_i$.  We can pick the largest such $i$, and get that $x_{i+1} \geq y_{i+1} > y_i$.  A move of the second type is to swap the portions of the two rows from column $i$ all the way to the far right edge of the board.  This results in an allowable arrangement as the only new pairs of adjacent numbers are that now $x_{i+1}$ is next to $y_i$ and $y_{i+1}$ is next to $x_i$, and by assumption, $x_{i+1} > y_i$ and $y_{i+1} > y_i > x_i$.  For the fourth condition, we have
\begin{eqnarray*}
r_i(x_i, x_{i-1}, \dots, x_g) & \leq & r_i(x_i, x_i-1, \dots,
x_i-i+1) = r_i(x_i+1)-1 \\ & < & r_i(x_i+1) \leq r_i(y_i) \leq
r_i(y_i, y_{i-1}, \dots, y_h).
\end{eqnarray*}

Suppose that $y_h > h > 1$.  Then for any $i < h$, a \textit{subdivision} is to replace the entries of the bottom row of columns $h$ through $i$ by $y_h-1, y_h-2, \dots, y_h-(h-i), y_h-(h-i)$.  This does not change the sums of conditions two, three, or four, as is easily seen by repeated application of the combinatorial identity ${n+1 \choose k} = {n \choose k} + {n \choose k-1}$.  Thus, if a subdivision is combined with other operations that satisfy conditions two, three, and four, and the end result satisfies condition one, it is an allowable move.

Suppose that at some point, the entries in columns $i$ through $g$ of the top row are $x_i, x_i-1, \dots, x_i-(i-g), x_i-(i-g)$.  Then a \textit{collapse} is to choose the largest value of $i$ with this property and replace these entries by an $x_i+1$ in column $i$ and clear all entries to the right of it. This is the inverse of a subdivision, and does not change the sums of conditions two, three, or four, so if combined with other operations satisfying those conditions in a way that ends with the first condition satisfied, it forms an allowable move.  Furthermore, the only new adjacent pair of entries that the collapse creates is $x_{i+1}$ adjacent to $x_i+1$.  Since $x_{i+1} > x_i+1$ by the choice of $i$, condition one will be satisfied in the top row for a move ending in a collapse.

Suppose that $y_h = h$ and $g = 1$.  A move of the third type starts by clearing the entry in column $h$ of the bottom row.  Increase $x_1$ by 1, and if necessary, do a collapse.  This operation strictly increases $r_{k+1}(x_k, \dots)$ without changing $r_{k+1}(y_k, \dots)$, so the sum of condition 3 strictly increases.  This is an allowable move because if the only changed entry not deleted were at least as large as the one to its left, it would have been fixed by a collapse.

Suppose that $y_h = h$ and $g > 1$.  A move of the fourth type is to clear the entry in column $h$ of the bottom row and put $g-1$ in column $g-1$ of the top row.  This is an allowable move because the only new entry is $g-1$ and immediately to its left is $x_g \geq g$.

Suppose that $x_h \geq y_h > h \geq g$ and $x_g - g < y_h - h$.  This implies $h > g$, as $h = g$ would yield $x_h \geq y_h > x_h$, a contradiction.  Pick the largest value of $i$ such that $y_h-(h-i) > x_i$; $i = g$ is such a value by assumption, so there must be a largest such value.  Since $x_h \geq y_h$, $i < h$.  By the choice of $i$, $x_{i+1} \geq y_h-(h-(i+1))$.  A move of the fifth type is to subdivide $y_h$ into $y_h-1, y_h-2, \dots, y_h-(h-i), y_h-(h-i)$, and then make a move of the second type to swap the two rows from column $i$ to the far right edge of the board.  Since $x_{i+1} \geq y_h-(h-(i+1)) > y_h-(h-i) > x_i$, we can make the move of the second type.  The new pairs of adjacent entries in the same row are $y_{h+1} > y_h-1 > y_h-2 > \dots > y_h-(h-i) > x_i$ in the bottom row and $x_{i+1} > y_h-(h-i)$ in the top row, so this is an allowable move.

Suppose that $y_h > h \geq g > 1$ and $x_g - g \geq y_h - h$.  The latter condition is equivalent to $x_g > y_h-(h-g+1)$.  A move of the sixth type is to subdivide $y_h$ into $y_h-1, y_h-2, \dots, y_h-(h-g+1), y_h-(h-g+1)$, and then move the $y_h-(h-g+1)$ from column $g-1$ of the bottom row to the top row. The new adjacent pairs are $x_g > y_h-(h-g+1)$ in the top row, and $y_{h+1} > y_h-1 > y_h-2 > \dots > y_h-(h-g+1)$ in the bottom row, so the first condition is satisfied.

Suppose that $g = h = 1$ and $x_1 \geq y_1$.  A move of the seventh type is to decrease $y_1$ by one (or delete it, if $y_1 = 1$), increase $x_1$ by 1, and collapse the top row as needed if $x_2 = x_1 + 1$.  Condition one is directly satisfied if $x_2 > x_1 + 1$. It is also satisfied if $x_2 = x_1 + 1$, as the move would end with a collapse. For condition three, we have
$${x_1 + 1 \choose 2} + {y_1 - 1 \choose 2} = {x_1 \choose 2} + {y_1 \choose 2} + x_1 - y_1 + 1 > {x_1 \choose 2} + {y_1 \choose 2},$$
so this is an allowable move, and the relevant sum strictly increases.

Suppose that $h > g = 1$ and $y_h - h \leq x_1 - 1$. The latter condition is equivalent to $y_h - h < x_1$.  A move of the eighth type is to subdivide the bottom row so that the entries from column $h$ on rightward become $y_h-1, y_h-2, \dots, y_h-(h-1), y_h-(h-1)$, and then make a move of the seventh type. The bottom row satisfies the first condition, as the new adjacent entries are $y_{h+1} > y_h-1 > y_h-2 > \dots > y_h-h$.  The top row also satisfies the first condition, as it is only changed by a move of the seventh type.  Hence, this is an allowable move.

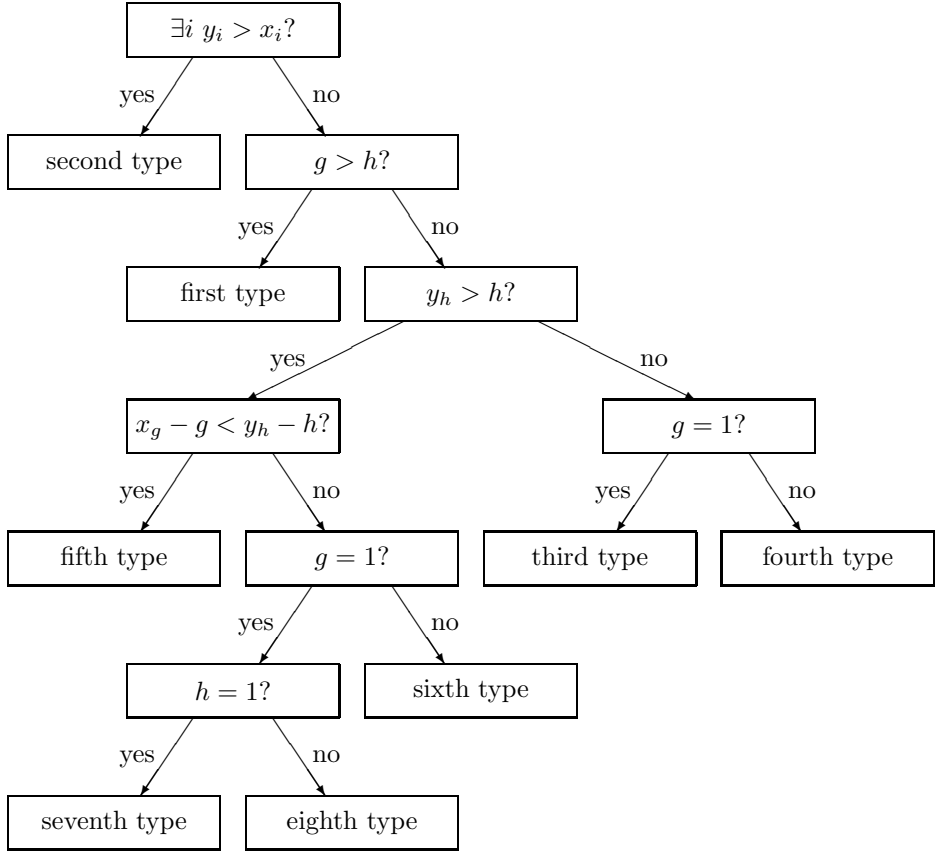
\begin{figure}[h]
\begin{picture}(350, 320)
\put (45, 300){\framebox(80, 20){$\exists i$ $y_i > x_i$?}}
\put (0, 250){\framebox(80, 20){second type}}
\put (90, 250){\framebox(80, 20){$g > h$?}}
\put (45, 200){\framebox(80, 20){first type}}
\put (135, 200){\framebox(80, 20){$y_h > h$?}}
\put (45, 150){\framebox(80, 20){$x_g - g < y_h - h$?}}
\put (225, 150){\framebox(80, 20){$g = 1$?}}
\put (0, 100){\framebox(80, 20){fifth type}}
\put (90, 100){\framebox(80, 20){$g = 1$?}}
\put (180, 100){\framebox(80, 20){third type}}
\put (270, 100){\framebox(80, 20){fourth type}}
\put (45, 50){\framebox(80, 20){$h = 1$?}}
\put (135, 50){\framebox(80, 20){sixth type}}
\put (0, 0){\framebox(80, 20){seventh type}}
\put (90, 0){\framebox(80, 20){eighth type}}
\put (70, 300){\vector(-2, -3){20}}
\put (100, 300){\vector(2, -3){20}}
\put (115, 250){\vector(-2, -3){20}}
\put (145, 250){\vector(2, -3){20}}
\put (150, 200){\vector(-2, -1){60}}
\put (200, 200){\vector(2, -1){60}}
\put (70, 150){\vector(-2, -3){20}}
\put (100, 150){\vector(2, -3){20}}
\put (250, 150){\vector(-2, -3){20}}
\put (280, 150){\vector(2, -3){20}}
\put (115, 100){\vector(-2, -3){20}}
\put (145, 100){\vector(2, -3){20}}
\put (70, 50){\vector(-2, -3){20}}
\put (100, 50){\vector(2, -3){20}}
\put (42, 283){yes}
\put (115, 283){no}
\put (87, 233){yes}
\put (160, 233){no}
\put (99, 183){yes}
\put (239, 183){no}
\put (42, 133){yes}
\put (115, 133){no}
\put (222, 133){yes}
\put (295, 133){no}
\put (87, 83){yes}
\put (160, 83){no}
\put (42, 33){yes}
\put (115, 33){no}
\end{picture}
\caption{Picking a type of move}
\end{figure}

If $x_k \leq a_k$, we can always make an allowable move, as shown in Figure~1. We then repeat the process and keep making such moves until $x_k > a_k$, at which point, we stop.  Condition 3 provides the weak inequality in the statement of the lemma.

The final move must increase $x_k$, and the only way that any of the eight types of moves used can change $x_k$ is to end with a collapse that goes all the way across the top row.  This ensures that the top row has only one entry.  The only three moves to potentially involve a collapse are those of the third, seventh, and eighth types, all of which strictly increase the sum of condition 3.  That the final move must be one of these types ensures that one of them is used in the algorithm, and so the inequality of the lemma is strict.

Further, a collapse only increases $x_k$ by one, so we now have $x_k = a_k+1$.  Condition two of allowable moves and the definition of $b_k, b_{k-1}, \dots$ then ensure that at this point, the only possible configuration of the board is for the entries of the bottom row to be $b_k, b_{k-1}, \dots$, which is what we wanted.

Condition one ensures that if $r_k(x_k, \dots, x_g) \geq r_k(a_k+1)$, then $x_k \geq a_k+1$, and the algorithm terminates.  Condition four says this sum must increase by at least one with each step.  Since the sum trivially cannot be negative, the algorithm then terminates in at most $r_k(a_k+1)$ steps. \endproof

Next is an easy extension of the previous lemma.

\begin{lemma} \label{algorithm2}
If $m = r_k(c_k, c_{k-1}, \dots) + r_k(a_k, a_{k-1}, \dots)$, $j > c_k$, $j > a_k$, and $m = r_k(j) + r_k(b_k, b_{k-1}, \dots)$, then
$$r_{k+1}(j) + r_{k+1}(b_k, \dots) > r_{k+1}(c_k, \dots) + r_{k+1}(a_k, \dots).$$
\end{lemma}

\proof  Assume without loss of generality that $c_k \geq a_k$. If we define $d_i$s such that
$$r_k(a_k, \dots) + r_k(c_k, \dots) = r_k(c_k+1) + r_k(d_k, d_{k-1}, \dots),$$
then Lemma~\ref{algorithm} states that
\begin{equation} \label{alg2eq}
r_{k+1}(a_k, \dots) + r_{k+1}(c_k, \dots) < r_{k+1}(c_k+1) + r_{k+1}(d_k, \dots).
\end{equation}

If $j = c_k+1$, then we are done.  Otherwise, $j > c_k+1$, so we repeat the process, using Lemma~\ref{algorithm} and increasing the $c_k+1$ term by 1 again, as many times as necessary to bring it up to $j$.  Since this operation increases the right hand side of (\ref{alg2eq}) each time, the result follows. \endproof

The following lemma is the key result in the proof of the main theorem.

\begin{lemma}  \label{lemma5}
Let $G$ be a graph with $\clique{k}{G} = m$ and let $m = r_k(n_k, n_{k-1}) + r_{k-1}(a_{k-1}, \dots, a_{k-s})$ be the unique representation of Lemma~\ref{lgbdlemma}.  If $G$ contains an $n_k$-clique, then $\clique{k+1}{G} \leq \lgbd{k}{m}.$
\end{lemma}

\proof  If $k < 3$, then $\lgbd{k}{m} = \oldbd{k}{m}$, so the lemma holds by the Kruskal-Katona theorem.  Otherwise, we can assume that $k \geq 3$.

Let $U$ be the vertex set of an $n_k$-clique of $G$, and $V$ the set of vertices of $G$ not contained in $U$. If the lemma is false, there must be a counterexample for which $|V|$ is minimal. If $|V| = 0$, then $G$ is an $n_k$-clique, which clearly satisfies the lemma.  If $|V| = 1$ and the one vertex of $V$ has degree $n_{k-1}$, then we have $\clique{k}{G} = r_k(n_k, n_{k-1})$ and $\clique{k+1}{G} = r_{k+1}(n_k, n_{k-1})$, which likewise satisfies the lemma.  Hence we must have $|V| \geq 2$.

Let $v \in V$.  Then $G - \{v\}$ has one fewer vertex in its own $V$ set, so it must satisfy the lemma.  A $k$-clique of $G$ either contains $v$ or it does not.  If it does, then it corresponds to the $(k-1)$-clique of \link{G}{v} consisting of the clique minus $v$.  If not, then it is a $k$-clique of $G - \{v\}$.  These correspondences reverse, so we have $\clique{k}{G} = \clique{k}{G - \{v\}} + \clique{k-1}{\link{G}{v}}$. By the same argument, $\clique{k+1}{G} = \clique{k+1}{G - \{v\}} + \clique{k}{\link{G}{v}}$.

Define $b_i$s and $c_i$s by
\begin{eqnarray*}
\clique{k}{G - \{v\}} & = & r_k(n_k, b_{k-1}, b_{k-2}, \dots) \spand \\ \clique{k-1}{\link{G}{v}} & = & r_{k-1}(c_{k-1}, c_{k-2}, \dots).
\end{eqnarray*}
The leading term of $\clique{k}{G - \{v\}}$ written in the form of Lemma~\ref{kklemma} is indeed $n_k$, as it contains a clique on $n_k$ vertices, and is a subgraph of $G$, so $\clique{k}{G - \{v\}} \leq \clique{k}{G} < r_k(n_k+1)$. By the Kruskal-Katona theorem,
\begin{eqnarray*}
\clique{k+1}{G - \{v\}} & \leq & r_{k+1}(n_k, b_{k-1}, \dots) \spand \\ \clique{k}{\link{G}{v}} & \leq & r_k(c_{k-1}, \dots).
\end{eqnarray*}
We are given that
\begin{eqnarray*}
\clique{k}{G} & = & r_k(n_k) + r_{k-1}(b_{k-1}, \dots) + r_{k-1}(c_{k-1}, \dots) \spand \\ \clique{k}{G} & = & r_k(n_k, n_{k-1}) + r_{k-1}(a_{k-1}, \dots).
\end{eqnarray*}

Applying the above inequalities gives
$$\clique{k+1}{G} = \clique{k+1}{G - \{v\}} + \clique{k}{\link{G}{v}} \leq r_{k+1}(n_k) + r_k(b_{k-1}, \dots) + r_k(c_{k-1}, \dots).$$
Then it suffices to show that
$$r_{k+1}(n_k) + r_k(b_{k-1}, \dots) + r_k(c_{k-1}, \dots) \leq r_{k+1}(n_k, n_{k-1}) + r_k(a_{k-1}, \dots),$$
or equivalently,
\begin{equation}\label{mainlem}
r_k(b_{k-1}, \dots) + r_k(c_{k-1}, \dots) \leq r_k(n_{k-1}) + r_k(a_{k-1}, \dots).
\end{equation}

Suppose that $n_{k-1} > b_{k-1}$ and $n_{k-1} > c_{k-1}$.  Then Lemma~\ref{algorithm2} immediately gives us (\ref{mainlem}).

Now suppose that $n_{k-1} = b_{k-1}$.  Since $G - \{v\}$ satisfies the lemma, we can define $d_i$s by
$$\clique{k}{G - \{v\}} = r_k(n_k, n_{k-1}) + r_{k-1}(d_{k-1}, d_{k-2}, \dots)$$
and have the bound
$$\clique{k+1}{G - \{v\}} \leq r_{k+1}(n_k, n_{k-1}) + r_k(d_{k-1}, d_{k-2}, \dots).$$
Then
\begin{eqnarray*}
\clique{k}{G} & = & \clique{k}{G - \{v\}} + \clique{k-1}{\link{G}{v}} \\ & = & r_k(n_k, n_{k-1}) + r_{k-1}(d_{k-1}, \dots) + r_{k-1}(c_{k-1}, \dots).
\end{eqnarray*}
Since
$$\clique{k}{G} = r_k(n_k, n_{k-1}) + r_{k-1}(a_{k-1}, \dots),$$
we obtain
$$r_{k-1}(a_{k-1}, \dots) = r_{k-1}(d_{k-1}, \dots) + r_{k-1}(c_{k-1}, \dots).$$
Hence by Lemma~\ref{disjoint},
$$r_k(a_{k-1}, \dots) \geq r_k(d_{k-1}, \dots) + r_k(c_{k-1}, \dots).$$
Putting the above inequalities together yields
\begin{eqnarray*}
\clique{k+1}{G} & = & \clique{k+1}{G - \{v\}} + \clique{k}{\link{G}{v}} \\ & \leq & r_{k+1}(n_k, n_{k-1}) + r_k(d_{k-1}, \dots) + r_k(c_{k-1}, \dots) \\ & \leq & r_{k+1}(n_k, n_{k-1}) + r_k(a_{k-1}, \dots) \\ & = & \lgbd{k}{m}.
\end{eqnarray*}

The remaining case is $n_{k-1} = c_{k-1} > b_{k-1}$.  We have not made any restrictions on the choice of $v$ except for $v \in V$, so if a different choice of $v$ puts us in one of the earlier cases, we are done.  That leaves only the case where $n_{k-1} = c_{k-1} > b_{k-1}$ regardless of the choice of $v$.

Pick vertices $p, q \in V$ and define graphs $H = G - \{p\}$ and $J = G - \{q\}$.  Assume without loss of generality that $\clique{k}{H} \geq \clique{k}{J}$.  Both $H$ and $J$ contain all $k$-cliques of $G$ that include neither $p$ nor $q$ as vertices.  By construction, the rest of the $k$-cliques of $H$ are those that contain $q$ but not $p$, and the remaining $k$-cliques of $J$ are those that contain $p$ but not $q$.  Since $\clique{k}{H} \geq \clique{k}{J}$, there are at least as many $k$-cliques of $G$ containing $q$ but not $p$ as vice versa.

Define $d_i$s and $e_i$s by
\begin{eqnarray*}
\clique{k}{H} & = & r_k(n_k, d_{k-1}, d_{k-2}, \dots) \spand \\ \clique{k-1}{\link{J}{p}} & = & r_{k-1}(e_{k-1}, e_{k-2}, \dots).
\end{eqnarray*}
Since $H$ contains all $k$-cliques of $G$ containing $q$ but not $p$ as well as all ${n_k \choose k}$ $k$-cliques of the $n_k$ vertices of $U$, and these are disjoint sets of cliques, the number of $k$-cliques of $G$ containing $q$ but not $p$ is at most $r_{k-1}(d_{k-1}, \dots)$. Each $(k-1)$-clique of \link{J}{p} corresponds to a $k$-clique of $G$ containing $p$ but not $q$, so there are $r_{k-1}(e_{k-1}, \dots)$ such cliques.  Thus,
\begin{equation}\label{deineq}
r_{k-1}(d_{k-1}, \dots) \geq r_{k-1}(e_{k-1}, \dots).
\end{equation}

If $p$ and $q$ are not connected by an edge, then by taking $v = p$, we get that $b_{k-1} \geq c_{k-1}$, a previous case. Otherwise, $p$ and $q$ must be connected by an edge, so we can define $f_i$s by
$$\clique{k-1}{\link{G}{pq}} = r_{k-1}(f_{k-2}, f_{k-3}, \dots).$$

A $(k-1)$-clique in the link of $pq$ has all vertices adjacent to $p$ in $G$, so if $p$ is added, it gives a $k$-clique containing $p$ but not $q$.  This is a $k$-clique in $J$ containing $p$, so it corresponds to a unique $(k-1)$-clique in \link{J}{p}.  Thus, $\clique{k-1}{\link{G}{pq}} \leq \clique{k-1}{\link{J}{p}}$, or equivalently,
\begin{equation}\label{efineq}
r_{k-1}(e_{k-1}, \dots) \geq r_{k-1}(f_{k-2}, \dots).
\end{equation}

By Lemma~\ref{dimshift},
$$r_{k-2}(e_{k-1}, \dots) \geq r_{k-2}(f_{k-2}, \dots).$$
Hence,
\begin{eqnarray}\label{eplus1}
r_{k-1}(e_{k-1}, \dots) + r_{k-2}(f_{k-2}, \dots) & \leq & r_{k-1}(e_{k-1}, \dots) + r_{k-2}(e_{k-1}, \dots) \nonumber \\ & = & r_{k-1}(e_{k-1}+1, e_{k-2}+1, \dots).
\end{eqnarray}
Applying Lemma~\ref{dimshift} to (\ref{deineq}) gives
$$r_{k-2}(d_{k-1}, \dots) \geq r_{k-2}(e_{k-1}, \dots).$$
Add the last inequality to (\ref{deineq}) to obtain
$$r_{k-1}(d_{k-1}, \dots) + r_{k-2}(d_{k-1}, \dots) \geq r_{k-1}(e_{k-1}, \dots) + r_{k-2}(e_{k-1}, \dots),$$
or equivalently,
\begin{equation}\label{dplus1}
r_{k-1}(d_{k-1}+1, d_{k-2}+1, \dots) \geq r_{k-1}(e_{k-1}+1, e_{k-2}+1, \dots).
\end{equation}

Let
\begin{eqnarray}
z & = & r_{k-1}(d_{k-1}, \dots) + r_{k-1}(e_{k-1}, \dots) + r_{k-2}(f_{k-2}, \dots) \label{zdef} \spand \\ z^+ & = & r_k(d_{k-1}, \dots) + r_k(e_{k-1}, \dots) + r_{k-1}(f_{k-2}, \dots).\label{z+def}
\end{eqnarray}
Applying the Kruskal-Katona theorem to the definitions of the $d_i$s, $e_i$s and $f_i$s yields
\begin{eqnarray*}
\clique{k+1}{H} & \leq & r_{k+1}(n_k, d_{k-1}, d_{k-2}, \dots), \\ \clique{k}{\link{J}{p}} & \leq & r_k(e_{k-1}, \dots), \spand \\ \clique{k-2}{\link{G}{pq}} & \geq & r_{k-2}(f_{k-2}, \dots).
\end{eqnarray*}
Applying these three inequalities to the definitions of $z$ and $z^+$ provides
\begin{eqnarray}
z + r_k(n_k) & \leq & \clique{k}{H} + \clique{k-1}{\link{J}{p}} + \clique{k-2}{\link{G}{pq}} = \clique{k}{G} \label{zbound} \ \textup{ and}\ \ \ \\ z^+ + r_{k+1}(n_k) & \geq & \clique{k+1}{H} + \clique{k}{\link{J}{p}} + \clique{k-1}{\link{G}{pq}} = \clique{k+1}{G}. \label{z+bound}
\end{eqnarray}

Suppose that $z \leq r_{k-1}(n_{k-1})$.  By (\ref{efineq}) and the proof of Lemma~\ref{linksub}, there is a complex with $r_{k-1}(e_{k-1}, \dots) + r_{k-2}(f_{k-2}, \dots)$ $(k-1)$-faces and $r_k(e_{k-1}, \dots) + r_{k-1}(f_{k-2}, \dots)$ $k$-faces. The disjoint union of this complex and the rev-lex complex on $r_k(d_{k-1}, \dots)$ $k$-faces has $z$ $(k-1)$-faces and $z^+$ $k$-faces. By the Kruskal-Katona theorem, since the complex has at most $r_{k-1}(n_{k-1})$ $(k-1)$-faces, it has at most $r_k(n_{k-1})$ $k$-faces.  Then $z^+ < r_k(n_{k-1})$, and so
$$\clique{k+1}{G} \leq z^+ + r_{k+1}(n_k) \leq r_{k+1}(n_k, n_{k-1}) \leq r_{k+1}(n_k, n_{k-1}) + r_k(a_{k-1}, \dots),$$
as desired.

Otherwise, $z > r_{k-1}(n_{k-1})$.  Define $g_i$s and $h_i$s such that
\begin{eqnarray}
z & = & r_{k-1}(n_{k-1}) + r_{k-1}(g_{k-1}, \dots) \label{zjg}  \spand \\ r_{k-1}(h_{k-1}, \dots) & = & r_{k-1}(e_{k-1}, \dots) + r_{k-2}(f_{k-2}, \dots). \nonumber
\end{eqnarray}
Substituting the latter into (\ref{zdef}) gives
\begin{equation}\label{zdh}
z = r_{k-1}(d_{k-1}, \dots) + r_{k-1}(h_{k-1}, \dots).
\end{equation}
By (\ref{efineq}), Lemma~\ref{linksub} gives
\begin{eqnarray}
& & r_{k+1}(n_k, d_{k-1}, \dots) + r_k(e_{k-1}, \dots) + r_{k-1}(f_{k-2}, \dots) \nonumber \\ & \leq & r_{k+1}(n_k, d_{k-1}, \dots) + r_k(h_{k-1}, \dots) \nonumber \\ & = & r_{k+1}(n_k, h_{k-1}, \dots) + r_k(d_{k-1}, \dots). \label{z+dh}
\end{eqnarray}

Recall that this case was based on the assumptions that $b_{k-1} < n_{k-1}$ and $c_{k-1} = n_{k-1}$.  Taking $v = p$ gives us $d_{k-1} = b_{k-1} < n_{k-1}$ and $h_{k-1} \leq c_{k-1} = n_{k-1}$.  If $h_{k-1} < n_{k-1}$, then combine (\ref{zjg}) and (\ref{zdh}) to get
$$r_{k-1}(d_{k-1}, \dots) + r_{k-1}(h_{k-1}, \dots) = r_{k-1}(n_{k-1}) + r_{k-1}(g_{k-1}, \dots).$$
Apply Lemma~\ref{algorithm2} and add $r_{k+1}(n_k)$ to both sides to obtain
\begin{equation}
r_{k+1}(n_k, h_{k-1}, \dots) + r_k(d_{k-1}, \dots) \leq r_{k+1}(n_k, n_{k-1}) + r_k(g_{k-1}, \dots). \label{nhdjg}
\end{equation}

Otherwise, $h_{k-1} = n_{k-1}$.  Since $d_{k-1} < n_{k-1}$, we must have $d_{k-1} + 1\leq n_{k-1}$.  Combine (\ref{eplus1}) and (\ref{dplus1}) to get
$$r_{k-1}(d_{k-1}+1, d_{k-2}+1, \dots) \geq r_{k-1}(e_{k-1},
\dots) + r_{k-2}(f_{k-2}, \dots) = r_{k-1}(h_{k-1}, \dots).$$
Then by Lemma~\ref{dimshift},
\begin{equation}
r_k(d_{k-1}+1, d_{k-2}+1, \dots) \geq r_k(h_{k-1}, \dots) \label{cond4}
\end{equation}
and $d_{k-1}+1 \geq h_{k-1}$.  This yields $d_{k-1}+1 \leq n_{k-1} = h_{k-1} \leq d_{k-1}+1$, so equality must hold throughout, and we have $d_{k-1} = n_{k-1} - 1$ and $h_{k-1} = n_{k-1}$.

Substitute these values of $d_{k-1}$ and $h_{k-1}$ into (\ref{zjg}) and (\ref{zdh}) and subtract $r_{k-1}(n_{k-1})$ to obtain
\begin{eqnarray*}
z - r_{k-1}(n_{k-1}) & = & r_{k-1}(g_{k-1}, g_{k-2}, \dots) \spand \\ z - r_{k-1}(n_{k-1}) & = & r_{k-2}(h_{k-2}, \dots) + r_{k-1}(n_{k-1}-1, d_{k-2}, \dots),
\end{eqnarray*}
so we have
\begin{equation}
r_{k-2}(h_{k-2}, \dots) + r_{k-1}(n_{k-1}-1, d_{k-2}, \dots) = r_{k-1}(g_{k-1}, g_{k-2}, \dots). \label{ghnd}
\end{equation}
Subtract $r_k(n_{k-1})$ from both sides of (\ref{cond4}) and use $n_{k-1}-1 > d_{k-2} > d_{k-3} > \dots$ to get
\begin{eqnarray*}
r_{k-1}(h_{k-2}, h_{k-3}, \dots) & \leq & r_{k-1}(d_{k-2} + 1, d_{k-3} + 1, \dots) \\ & \leq & r_{k-1}(n_{k-1}-1, d_{k-2}, \dots).
\end{eqnarray*}
Applying Lemma~\ref{linksub} to (\ref{ghnd}) yields
$$r_{k-1}(h_{k-2}, \dots) + r_k(n_{k-1}-1, d_{k-2}, \dots) \leq r_k(g_{k-1}, g_{k-2}, \dots).$$
Add $r_{k+1}(n_k, n_{k-1})$ to both sides to get (\ref{nhdjg}) in this case also.

Finally, we chain together (\ref{z+bound}), (\ref{z+def}), (\ref{z+dh}), (\ref{nhdjg}), the definition of \lgbd{k}{m}, (\ref{zjg}), and (\ref{zbound}) to conclude
\begin{eqnarray*}
\clique{k+1}{G} & \leq & z^+ + r_{k+1}(n_k) \\ & = & r_{k+1}(n_k, d_{k-1}, \dots) + r_k(e_{k-1}, \dots) + r_{k-1}(f_{k-2}, \dots) \\ & \leq & r_{k+1}(n_k, h_{k-1}, \dots) + r_k(d_{k-1}, \dots) \\ & \leq & r_{k+1}(n_k, n_{k-1}) + r_k(g_{k-1}, \dots) \\ & = & \lgbd{k}{r_k(n_k, n_{k-1}) + r_{k-1}(g_{k-1}, \dots)} \\ & = & \lgbd{k}{r_k(n_k) + z} \\ \hspace{76 pt} & \leq & \lgbd{k}{\clique{k}{G}}. \hspace{198 pt} \square
\end{eqnarray*}

We are now ready to prove our main theorem, Theorem~\ref{maintheorem}, which states that for a graph $G$ with $\clique{k}{G} = m$,
$$\clique{k+1}{G} \leq \max\{\lgbd{k}{m}, \smbd{k}{m}\}.$$

\proof  Let $\clique{k}{G} = r_k(n_k, n_{k-1}, \dots)$.  Either $G$ has a clique on $n_k$ vertices or it does not.  If $G$ does have a clique on $n_k$ vertices, then by Lemma~\ref{lemma5},
$\clique{k+1}{G} \leq \lgbd{k}{m} \leq \max\{\lgbd{k}{m}, \smbd{k}{m}\}$.  If $G$ does not have a clique on $n_k$ vertices, then by Theorem~\ref{mytheorem}, $\clique{k+1}{G} \leq \smbd{k}{m} \leq \max\{\lgbd{k}{m}, \smbd{k}{m}\}$. \endproof

\section {Evaluating the bound}

In this section, we evaluate the bound of Theorem~\ref{maintheorem}.  There are several questions to consider.  First, which of the two bounds is larger?  Next, how close is each bound to being sharp?  Finally, how close is each bound to the bound of the Kruskal-Katona theorem?

On the first question, we typically have $\lgbd{k}{m} > \smbd{k}{m}$ if $n_k$ is close to $k$, if $a_{k-1}$ is small relative to $n_k$ (which must happen if $n_{k-1}$ is small relative to $n_k$), or if $n_{k-1} = n_k-1$, and $\lgbd{k}{m} < \smbd{k}{m}$ otherwise. Empirically, for fixed $k$ with $k$ small, if we define a sequence $f_j = {\#\{m \leq j\ |\ \lgbd{k}{m} > \smbd{k}{m}\} \over j}$, this sequence seems to converge to a number around 0.7.  That is, $\lgbd{k}{m} > \smbd{k}{m}$ a substantial majority of the time.  Because there can be many consecutive values of $m$ for which $n_{k-1}$ is much smaller than $n_k$, there are very long sets of consecutive increasing terms of the sequence $f_j$, making the sequence converge very slowly, if it converges at all.

Regarding the question of how good each of the bounds is, we evaluate the \lgbd{k}{m} and \smbd{k}{m} bound separately.  Loosely, the \lgbd{k}{m} bound is not sharp, but is much better than Kruskal-Katona.  The \smbd{k}{m} bound is also an improvement over Kruskal-Katona, but not necessarily such a good bound when it applies.  We also give several constructions of graphs that attain the new bound under certain circumstances.

It is worth mentioning that Theorem~\ref{maintheorem} characterizes exactly when the bounds of the Kruskal-Katona theorem are attained by a graph for $k \leq 7$.  In this case, the bound is attained by a graph exactly when $\lgbd{k}{m} = \oldbd{k}{m}$.  The proof of this fact is a long computation that we decided to omit.

\subsection {The $\lgbd{k}{m}$ bound}

In this section, we evaluate the bound of Lemma~\ref{lemma5}.  The inequality $\lgbd{k}{m} \leq \oldbd{k}{m}$ follows immediately from Lemma~\ref{disjoint}.  Still, that leaves open the question of whether it is much less, or whether the inequality is even strict. Theorem~\ref{convprob} is a convergence in probability type of result that shows that \lgbd{k}{m} is almost always much closer to being sharp than it is to \oldbd{k}{m}.

Throughout this section, let $m = r_k(n_k, n_{k-1}, \dots)$ be the representation of Lemma~\ref{kklemma} and let $m = r_k(n_k, n_{k-1}) + r_{k-1}(a_{k-1}, \dots)$ be the representation of Lemma~\ref{lgbdlemma}.   First we need a couple definitions in order to state the result more precisely.

\begin{definition}
\textup{Let \conbd{k}{m} be the largest number of $(k+1)$-cliques that a graph with $m$ $k$-cliques and at least one $n_k$-clique can possibly have.}
\end{definition}

While we do not have a formula for \conbd{k}{m}, Lemma~\ref{lemma5} states that $\conbd{k}{m} \leq \lgbd{k}{m}$.  We can also readily get a constructive lower bound.

\begin{lemma}
For all $m, k > 0$, $\conbd{k}{m} \geq r_{k+1}(n_k, n_{k-1}) + r_k(a_{k-1}).$
\end{lemma}

\proof Start with a clique on $n_k$ vertices.  Add a new vertex and make it adjacent to $n_{k-1}$ previous vertices.  Add another new vertex and make it adjacent to $a_{k-1}$ of the first $n_k$ vertices.  If we call this graph $G$, then we have
\begin{eqnarray*}
\clique{k}{G} & = & r_k(n_k, n_{k-1}) + r_{k-1}(a_{k-1}) \leq m \spand \\ \hspace{60 pt} \clique{k+1}{G} & = & r_{k+1}(n_k, n_{k-1}) + r_k(a_{k-1}) \leq \conbd{k}{m}. \hspace{60 pt} \square
\end{eqnarray*}

\begin{definition}
\textup{Given $m$ and $k$ with $\conbd{k}{m} \not = \oldbd{k}{m}$, define
$$\ratio{k}{m} := {\lgbd{k}{m} - \conbd{k}{m} \over \oldbd{k}{m} - \conbd{k}{m}}.$$
If $\conbd{k}{m} = \oldbd{k}{m}$, define $\ratio{k}{m} := 1$.}
\end{definition}

It is immediate from the definition that $0 \leq \ratio{k}{m} \leq 1$.  Intuitively \lgbd{k}{m} is a good bound when \ratio{k}{m} is small.  When $\conbd{k}{m} = \oldbd{k}{m}$, the \lgbd{k}{m} bound is irrelevant, so choice of \ratio{k}{m} is arbitrary and does not affect the final result.

First we need a lemma that is a very rough approximation.

\begin{lemma} \label{ratbound}
If $n_{k-2} > k^2$, then $\ratio{k}{m} \leq {k^2 \over n_{k-2} - k^2}$.
\end{lemma}

\proof  We can compute
\begin{equation}
r_{k-2}(n_{k-2}+1) > r_{k-2}(n_{k-2}, \dots) = m - r_k(n_k, n_{k-1}) = r_{k-1}(a_{k-1}, \dots).\label{naineq}
\end{equation}
We must have $k \geq 3$ in order for $n_{k-2}$ to be defined, so $n_{k-2} > k^2 > 2k$.  If $a_{k-1} > n_{k-2}$, then $n_{k-2} > 2k$ and (\ref{naineq}) would yield
$$r_{k-1}(a_{k-1}) \geq r_{k-1}(n_{k-2}+1) > r_{k-2}(n_{k-2}+1) > r_{k-1}(a_{k-1}, \dots),$$
a contradiction.  Hence, $a_{k-1} \leq n_{k-2}$.

We can also use (\ref{naineq}) to get
\begin{eqnarray*}
{n_{k-2}+1 \choose k-2} & > & {a_{k-1} \choose k-1}, \qquad \mbox{and hence} \\ {n_{k-2} \choose k-1} & > & {k(n_{k-2}-k+2)(n_{k-2}-k+3) \over (k-1)(a_{k-1}-k+1)(n_{k-2}+1)}{a_{k-1} \choose k}.
\end{eqnarray*}
We can then apply the definition of \oldbd{k}{m} to obtain
\begin{eqnarray*}
\oldbd{k}{m} & \geq & r_{k+1}(n_k, n_{k-1}, n_{k-2}) \\ & > & r_{k+1}(n_k, n_{k-1}) + {k(n_{k-2}-k+2)(n_{k-2}-k+3) \over (k-1)(a_{k-1}-k+1)(n_{k-2}+1)}\ r_k(a_{k-1}).
\end{eqnarray*}

The quantity \ratio{k}{m} is made larger if we overestimate \lgbd{k}{m} and if we underestimate \oldbd{k}{m} and \conbd{k}{m}. Then we can use the bound of the previous paragraph as well as
\begin{eqnarray*}
\lgbd{k}{m} & \leq & r_{k+1}(n_k, n_{k-1}) + r_k(a_{k-1}+1) \spand \\ \conbd{k}{m} & \geq & r_{k+1}(n_k, n_{k-1}) + r_k(a_{k-1})
\end{eqnarray*}
to get an upper bound on \ratio{k}{m} of
$${r_{k+1}(n_k, n_{k-1}) + r_k(a_{k-1}+1) - (r_{k+1}(n_k, n_{k-1}) + r_k(a_{k-1})) \over r_{k+1}(n_k, n_{k-1}) + {k(n_{k-2}-k+2)(n_{k-2}-k+3) \over (k-1)(a_{k-1}-k+1)(n_{k-2}+1)}r_k(a_{k-1}) - (r_{k+1}(n_k, n_{k-1}) + r_k(a_{k-1}))} $$
\begin{eqnarray*}
& = & {r_{k-1}(a_{k-1}) \over {k(n_{k-2}-k+2)(n_{k-2}-k+3) \over (k-1)(a_{k-1}-k+1)(n_{k-2}+1)}r_k(a_{k-1}) - r_k(a_{k-1})} \\ & = & {{k \over a_{k-1}-k+1}r_k(a_{k-1}) \over \Big({k(n_{k-2}-k+2)(n_{k-2}-k+3) \over (k-1)(a_{k-1}-k+1)(n_{k-2}+1)}-1\Big)r_k(a_{k-1})} \\ & = & {{k \over a_{k-1}-k+1}(k-1)(a_{k-1}-k+1)(n_{k-2}+1) \over k(n_{k-2}-k+2)(n_{k-2}-k+3) - (k-1)(a_{k-1}-k+1)(n_{k-2}+1)} \\ & \leq & {k(k-1)(n_{k-2}+1) \over k(n_{k-2}-k+2)(n_{k-2}-k+3) - (k-1)(n_{k-2}-k+1)(n_{k-2}+1)} \\ & = & {k(k-1)(n_{k-2}+1) \over (n_{k-2}+1)(n_{k-2}-k^2) + k^2(k-3) + 2kn_{k-2} + n_{k-2} + 4k + 1} \\ & \leq & {k(k-1)(n_{k-2}+1) \over (n_{k-2}+1)(n_{k-2}-k^2)} \leq  {k^2 \over n_{k-2}-k^2}. \hspace{160 pt} \square
\end{eqnarray*}

\begin{lemma} \label{ak2bound}
For every $k \geq 3$ and $w$,
$$\lim_{j \to \infty} {\#\{m \leq j\ |\ n_{k-2} < w\} \over j} = 0.$$
\end{lemma}

\proof  For any particular values of $n_k$ and $n_{k-1}$, there are at most ${w \choose k-2}$ corresponding values of $m$ with $n_{k-2} < w$.  If we define $n$ such that ${n \choose k} \leq j < {n+1 \choose k}$, then there are at most ${n \choose 2}$ ways to pick $n_k$ and $n_{k-1}$ corresponding to some value of $m \leq j$.  Hence, $\#\{m \leq j\ |\ n_{k-2} < w\} \leq {n \choose 2}{w \choose k-2}$.  Since $j \geq {n \choose k}$, we have
$${\#\{m \leq j\ |\ n_{k-2} < w\} \over j} \leq {{n \choose 2}{w \choose k-2} \over {n \choose k}}.$$
The right hand side is a rational function in $n$, with the numerator of degree 2 and the denominator of degree $k \geq 3$, so it goes to zero as $n \to \infty$.  If we let $j \to \infty$, then $n \to \infty$ as well, so we have
$$\limsup_{j \to \infty} {\#\{m \leq j\ |\ n_{k-2} < w\} \over j} \leq \limsup_{n \to \infty} {{n \choose 2}{w \choose k-2} \over {n \choose k}} = 0.$$
Since the lim sup of this non-negative sequence is not positive, the sequence must converge to zero. \endproof

Finally we reach the main result of this section.  With suitable definitions of distributions, it essentially says that \ratio{k}{m} converges to zero in probability.

\begin{theorem} \label{convprob}
For every $k \geq 3$ and $\epsilon > 0$,
$$\lim_{j \to \infty} {\#\{m \leq j\ |\ \ratio{k}{m} > \epsilon\} \over j} = 0.$$
\end{theorem}

\proof  Let $w = \Big\lceil {k^2(1+\epsilon) \over \epsilon} \Big\rceil$, so that ${k^2 \over w - k^2} \leq \epsilon$.  By Lemma~\ref{ratbound}, if $n_{k-2} \geq w$, then
$$\ratio{k}{m} \leq {k^2 \over n_{k-2} - k^2} \leq {k^2 \over w - k^2} \leq \epsilon.$$
Then by Lemma~\ref{ak2bound},
$$\hspace{23 pt} \limsup_{j \to \infty} {\#\{m \leq j\ |\ \ratio{k}{m} > \epsilon\} \over j} \leq \limsup_{j \to \infty} {\#\{m \leq j\ |\ n_{k-2} < w\} \over j} = 0.  \hspace{23 pt} \square$$

\subsection {The $\smbd{k}{m}$ bound}

In this section, we evaluate the \smbd{k}{m} bound.  Unlike the case of the \lgbd{k}{m} bound, if we were to define something analogous to \ratio{k}{m} here, it does not empirically seem to converge in probability.  It may converge weakly to some distribution, but this would be difficult to calculate,
and if we restrict to values of $m$ such that $\smbd{k}{m} > \lgbd{k}{m}$, it may not still converge to the same distribution.

Instead, we prove that $\smbd{k}{m} < \oldbd{k}{m}$ whenever $\oldbd{k}{m} > 0$, or equivalently, whenever \smbd{k}{m} is defined.  Hence, the (non-zero) bounds of the Kruskal-Katona theorem are never attained by a graph lacking the largest clique it could possibly have for its prescribed number of cliques of a given size.  Whether the bound of Theorem~\ref{maintheorem} is strictly tighter than that of the Kruskal-Katona theorem then depends only on the \lgbd{k}{m} bound.

First we need a lemma showing that the inequality of Lemma~\ref{linksub} is strict if we strengthen one assumption.

\begin{lemma} \label{strictlink}
If $r_k(c_k, c_{k-1}, \dots) = r_k(a_k, a_{k-1}, \dots) + r_{k-1}(b_{k-1}, b_{k-2}, \dots)$ and $c_k-1 = a_k > b_{k-1}$, then
$$r_{k+1}(c_k, \dots) > r_{k+1}(a_k, \dots) + r_k(b_{k-1}, \dots).$$
\end{lemma}

\proof  Subtract $r_k(a_k)$ from both sides of the equation of the lemma to get
$$r_{k-1}(a_k) + r_{k-1}(c_{k-1}, \dots) = r_{k-1}(a_{k-1}, \dots) + r_{k-1}(b_{k-1}, \dots).$$
Lemma~\ref{algorithm2} states that
$$r_k(a_k) + r_k(c_{k-1}, \dots) > r_k(a_{k-1}, \dots) + r_k(b_{k-1}, \dots).$$
Adding $r_{k+1}(a_k)$ to both sides completes the proof.
\endproof

\begin{proposition} \label{smallstrict}
If $m$ and $k$ satisfy $\oldbd{k}{m} > 0$, then $\smbd{k}{m} < \oldbd{k}{m}$.
\end{proposition}

\proof  Let $m = r_k(n_k, n_{k-1}, \dots)$ and $m = {a_k\choose k}_{n_k-1} + {a_{k-1}\choose k-1}_{n_k-2} + \dots$ be the representations used in the definitions of $\oldbd{k}{m}$ and $\smbd{k}{m}$.  Since ${n \choose k}_r \leq {n \choose k}$, $m \leq r_k(a_k, a_{k-1}, \dots)$, and so $a_k \geq n_k$.

Suppose that $n_k = a_k$.  Since $a_k - \big\lfloor{a_k\over a_k-1}\big\rfloor > a_{k-1}$ and $\big\lfloor{a_k\over a_k-1}\big\rfloor \geq 1$, we have that $n_k - 1 = a_k - 1 > a_{k-1}$, and so $n_k-2 \geq a_{k-1}$. Then we have that ${a_{k-1} \choose k-1}_{n_k-2} = {a_{k-1} \choose k-1}$, and similarly for all $a_{k-i}$ terms with $i \geq 1$.

The Tur\'{a}n graph $T_{a_k,a_k-1}$ consists of a clique on $a_k-1$ vertices and one other vertex adjacent to $a_k-2$ vertices.  Thus, ${a_k\choose k}_{a_k-1} = r_k(a_k-1,a_k-2)$, and so $m = r_k(a_k-1,a_k-2) + r_{k-1}(a_{k-1}, a_{k-2}, \dots)$.  Since $a_k-1 > a_{k-1}$, we get the desired inequality from Lemma~\ref{strictlink}.

Now we proceed by induction on $a_k$.  Let $C$ be the $r$-colored rev-lex complex on $m$ $k$-faces and \smbd{k}{m} $(k+1)$-faces and let $v$ be the first vertex in the rev-lex order.  The number of $k$-faces of $C$ is the number of those containing $v$ plus the number not containing $v$, so
\begin{eqnarray*}
m & = & \clique{k}{C} = \clique{k}{C - \{v\}} + \clique{k-1}{\link{C}{v}} \spand \\ \smbd{k}{m} & = & \clique{k+1}{C} = \clique{k+1}{C - \{v\}} + \clique{k}{\link{C}{v}}.
\end{eqnarray*}
Define $b_i$s and $c_i$s by
\begin{eqnarray*}
r_k(b_k, b_{k-1}, \dots) & = & \clique{k}{C - \{v\}} \spand \\ r_k(c_k, c_{k-1}, \dots) & = & \clique{k}{\link{C}{v}}.
\end{eqnarray*}
By construction, $C - \{v\}$ contains the clique complex of a $T_{a_k-1,r}$ and is contained in the clique complex of a $T_{a_k,r}$ but does not have a $k$-face on the last $k$ vertices.  Hence, if \clique{k}{C - \{v\}} is written as in the definition of \smbd{k}{\clique{k}{C - \{v\}}}, its leading term will be $a_k-1$.  As such, the inductive hypothesis provides
$$r_{k+1}(b_k, b_{k-1}, \dots) > \clique{k+1}{C - \{v\}}.$$
Since $\link{C}{v} \subset (C - \{v\})$, $r_k(b_k, \dots) \geq r_k(c_k, \dots)$, so by Lemma~\ref{linksub}, we have
\begin{eqnarray*}
\oldbd{k}{m} & \geq & r_{k+1}(b_k, \dots) + r_k(c_k, \dots) \\ & > & \clique{k+1}{C - \{v\}} + \clique{k}{\link{C}{v}} \\ \hspace{119pt} & = & \smbd{k}{m}. \hspace{169pt} \square
\end{eqnarray*}

\subsection {Attaining the bounds}

In this section, we give a few conditions under which the bounds of Theorem~\ref{maintheorem} are attained.  The proofs that the bounds are attained are by construction.

\begin{construction} \label{const1}
\textup{Let $m = r_k(n_k, n_{k-1}) + r_{k-1}(a_{k-1}, a_{k-2}, \dots)$ be the representation of $m$ satisfying the conditions of Lemma~\ref{lgbdlemma}.  Suppose that $a_{k-2} = k-2$ or does not exist.  Construct $G$ by starting with a clique on $n_k$ vertices. Add a new vertex and make it adjacent to $n_{k-1}$ of the original $n_k$ vertices.  Add another new vertex and make it adjacent to $a_{k-1}$ of the original vertices.  Then $\clique{k+1}{G} = \lgbd{k}{m}$ and $\clique{k}{G} \leq m$.  If the latter inequality is strict, we can add however many isolated cliques on $k$ vertices are needed to make equality hold.}
\end{construction}

\begin{construction} \label{const2}
\textup{Let $m = r_k(n_k, n_{k-1}) + r_{k-1}(a_{k-1}, a_{k-2}, \dots)$ be the representation of $m$ satisfying the conditions of Lemma~\ref{lgbdlemma}.  Suppose that $a_{k-3} = k-3$ or does not exist and $n_k + a_{k-2} \geq n_{k-1} + a_{k-1}$.  Construct $G$ by starting with a clique on $n$ vertices.  Add a new vertex $v$ and make it adjacent to $n_{k-1}$ of the original $n_k$ vertices.  Add another new vertex $u$ and make it adjacent to $a_{k-1}$ of the original vertices as well as adjacent to $v$ such that $u$ and $v$ are adjacent to $a_{k-2}$ common vertices.}

\textup{This can be done if $u$ is adjacent to $a_{k-2}$ vertices to which $v$ is also adjacent, and $a_{k-1} - a_{k-2}$ vertices (other than $v$) to which $v$ is not adjacent.  Since $n_{k-1} \geq a_{k-1} > a_{k-2}$, we can make the last two vertices adjacent to enough common neighbors.  We can prevent them from being adjacent to too many common neighbors if there are at least $n_{k-1} + a_{k-1} - a_{k-2}$ vertices in the first $n$ available.  That is, this construction can be done if $n_k \geq n_{k-1} + a_{k-1} - a_{k-2}$, or equivalently, $n_k + a_{k-2} \geq n_{k-1} + a_{k-1}$, the condition of the lemma.}

\textup{Then $\clique{k+1}{G} = \lgbd{k}{m}$ and $\clique{k}{G} \leq m$.  If the latter inequality is strict, we can add several isolated cliques on $k$ vertices to make equality hold.}
\end{construction}

\begin{construction}
\textup{Let $m = {a_k\choose k}_{n_k-1} + {a_{k-1}\choose k-1}_{n_k-2} + \dots + {a_{k-s}\choose k-s}_{n_k-s-1}$ be the representation satisfying the conditions of Lemma~\ref{colorcan}.  Suppose that $a_{k-2} = k-2$ or does not exist.  Let $G$ be the Tur\'{a}n graph $T_{a_k,n_k-1}$.  If we remove a part tied for the smallest from $G$, it still has at least $a_k-\big\lfloor{a_k\over n_k-1}\big\rfloor$ vertices remaining. Since $a_k-\big\lfloor{a_k \over n_k-1}\big\rfloor > a_{k-1}$,  $T_{a_k,n_k-1}$ has a Tur\'{a}n graph $T_{a_{k-1},n_k-2}$ as an induced subgraph.  Hence, we can create a graph $G'$ from $G$ by adding a new vertex adjacent to the vertices of a $T_{a_{k-1},n_k-2}$ induced subgraph of $G$. Then $\clique{k}{G'} \leq m$ and $\clique{k+1}{G'} = \smbd{k}{m}$. If the inequality is strict, we can add some isolated cliques on $k$ vertices to make equality hold.}
\end{construction}

\section {Bound for non-consecutive dimensions}

\begin{theorem} \label{nonconsec}
Let $k, i > 0$, $\clique{k}{G} = r_k(n_k, n_{k-1}) + r_{k-1}(a_{k-1}, a_{k-2}, \dots)$ be the representation of Lemma~\ref{lgbdlemma}, and $\clique{k}{G} = {b_k\choose k}_{n_k-1} + {b_{k-1}\choose k-1}_{n_k-2} + \dots + {b_{k-s}\choose k-s}_{n_k-s-1}$ be the representation of Lemma~\ref{colorcan}.  Then
$$c_{k+i}(G) \leq \max \bigg\{
\begin{array}{l}
r_{k+i}(n_k, n_{k-1}) + r_{k+i-1}(a_{k-1}, \dots) \\ {b_k\choose k+i}_{n_k-1} + \dots + {b_{k-s}\choose k+i-s}_{n_k-s-1}.
\end{array}$$
\end{theorem}

\proof  If $G$ has a clique on $n_k$ vertices, then by Lemma~\ref{lemma5}, $\clique{k+1}{G} \leq \lgbd{k}{m}$, $\clique{k+2}{G} \leq \lgbd{k+1}{\clique{k+1}{G}} \leq \lgbd{k+1}{\lgbd{k}{\clique{k}{G}}}$, and so forth, until we get
\begin{eqnarray*}
\clique{k+i}{G} & \leq & \lgbd{k+i-1}{\dots \lgbd{k+1}{\lgbd{k}{m}}\dots} \\ & = & r_{k+i}(n_k, n_{k-1}) + r_{k+i-1}(a_{k-1}, a_{k-2}, \dots).
\end{eqnarray*}

Otherwise, $G$ does not have a clique on $n_k$ vertices, in which case, by Theorem~\ref{mytheorem}, $\clique{k+i}{G} \leq {b_k\choose k+i}_{n_k-1} + \dots + {b_{k-s}\choose k+i-s}_{n_k-s-1}$.  Either way, the assertion holds.  \endproof

The next example shows that if $m = \clique{k}{G}$, then the first bound in Theorem~\ref{nonconsec} is \lgbd{k+i-1}{\dots \lgbd{k+1}{\lgbd{k}{m}}\dots}.  The second bound, on the other hand, may not be \smbd{k+i-1}{\dots \smbd{k+1}{\smbd{k}{m}}\dots}, but Theorem~\ref{nonconsec} could have been stated using this quantity as the second bound.

\begin{example}
\textup{Let $\clique{3}{G} = r_3(6) = 20$.  We can compute $20 = {6 \choose 3}_5 + {3 \choose 2}_4 + {1 \choose 1}_3$.  Then the second possible upper bound on \clique{5}{G} in Theorem~\ref{nonconsec} is ${6 \choose 5}_5 + {3 \choose 4}_4 + {1 \choose 3}_3 = 2$, while $\smbd{3}{20} = {6 \choose 4}_5 + {3 \choose 3}_4 + {1 \choose 2}_3 = 10$.  We can also compute $10 = r_4(5, 4) + r_3(3)$ and then $10 = {7 \choose 4}_4 + {4 \choose 3}_3$.  Then we compute $\smbd{4}{10} = {7 \choose 5}_4 + {4 \choose 4}_3 = 0$, so that $\smbd{4}{\smbd{3}{20}} = 0 \not = 2$.}

\textup{What happened here is that the value of $n_k$ used for the first time we apply the bound was not the same as for the second. This cannot be used to get an improved bound, though, as when this happens, the latter bound is always $< r_{k+i}(n_k)$, while the former bound is always $\geq r_{k+i}(n_k)$, and hence larger. By this logic, Theorem~\ref{nonconsec} could have instead stated that
$$c_{k+i}(G) \leq \max\{\lgbd{k+i-1}{\dots \lgbd{k+1}{\lgbd{k}{m}}\dots}, \smbd{k+i-1}{\dots
\smbd{k+1}{\smbd{k}{m}}\dots}\},$$
as this would never change the larger of the two bounds.}
\end{example}

Theorem~\ref{nonconsec} does sometimes give us a sharper bound than Theorem~\ref{maintheorem} alone, as shown in the next example.

\begin{example} \label{example}
\textup{If $\clique{3}{G} = 70$, then the bound of Theorem~\ref{maintheorem} is $\clique{4}{G} \leq \max\{\smbd{3}{70}, \lgbd{3}{70}\} = \max\{85, 81\} = 85$. This bound is attained by the Tur\'{a}n graph $T_{9,7}$, as $\clique{3}{T_{9,7}} = 70$ and $\clique{4}{T_{9,7}} = 85$.}

\textup{If $\clique{4}{G} = 85$, then by Theorem~\ref{maintheorem}, $\clique{5}{G} \leq  \max\{\smbd{4}{85}, \lgbd{4}{85}\} = \max\{61, 62\} = 62$. This bound is attained by Construction~\ref{const1}.}

\textup{By Theorem~\ref{nonconsec}, if $\clique{3}{G} = 70$, then $\clique{5}{G} \leq \max\{61, 61\} = 61$. Hence, there is a graph $G_1$ with $\clique{3}{G_1} = 70$ and $\clique{4}{G_1} = 85$ and there is a graph $G_2$ with $\clique{4}{G_2} = 85$ and $\clique{5}{G_2} = 62$, but there is no graph $G_3$ with $\clique{3}{G_3} = 70$ and $\clique{5}{G_3} = 62$.}
\end{example}

\textit{Acknowledgements.} I would like to thank my thesis advisor Isabella Novik for her many useful discussions and helpful advice throughout the writing and editing of this article.


\begin{thebibliography}{99}
\bibitem{colorfixed} L. Billera and A. Bjorner, Face numbers of polytopes and complexes, Handbook of Discrete and Computational Geometry, J.E.~Goodman and J.~O'Rourke, eds., CRC Press, Boca Raton, New York, 1997, pp. 291-310.
\bibitem{mainconj} J. Eckhoff, Intersection properties of boxes. I. An upper-bound theorem, Israel J. Math. 62 (1988), 283-301.
\bibitem{eckhoff0} J. Eckhoff, The maximum number of triangles in a $K_4$-free graph, Discrete Math. 194 (1999), 95-106.
\bibitem{eckhoff} J. Eckhoff, A new Tur\'{a}n-type theorem for cliques in graphs, Discrete Math. 282 (2004), 113-122.
\bibitem{balanced} P. Frankl, Z. F\"{u}redi, and G. Kalai, Shadows of colored complexes, Math. Scand. 63 (1988), 169-178.
\bibitem{koszul} R. Fr\"{o}berg, Koszul Algebras, in: Advances in Commutative Ring Theory, Lecture Notes in Pure and Appl. Math., Dekker, New York, 1999, pp. 337-350.
\bibitem{previous} A. Frohmader, Face vectors of flag complexes, to appear in Israel J. Math, math.CO/0605673
\bibitem{hhmtz} J. Herzog, T. Hibi, S. Murai, N. Trung, and X. Zhang, Kruskal-Katona type theorems for clique complexes arising from chordal and strongly chordal graphs, preprint, math.CO/0606477
\bibitem{katona} G. Katona, A theorem of finite sets, in: Theory of Graphs, Academic Press, New York, 1968, pp. 187-207.
\bibitem{kruskal} J.B. Kruskal, The number of simplices in a complex, in: Mathematical Optimization Techniques, University of California Press, Berkeley, California, 1963, pp. 251-278.
\bibitem{lovasz} L. Lov\'{a}sz, M. Simonovits,  On the number of complete subgraphs of a graph. II. Studies in pure  mathematics, 459-495, Birkhuser, Basel, 1983.
\bibitem{cohenface} R. Stanley, Cohen-Macaulay complexes, in: (M. Aigner, ed.) Higher Combinatorics, Reidel, Dordrecht, 1977.
\bibitem{greenbook} R. Stanley, Combinatorics and Commutative Algebra, Second Edition, Birkhauser Boston, Inc., Boston, Massachusetts, 1996, 53-64.
\bibitem{turan} P. Tur\'{a}n, Eine Extremalaufgabe aus der Graphentheorie, Mat. Fiz. Lapok 48 (1941), 436-452.
\bibitem{mantel} J. van Lint and R. Wilson, A Course in Combinatorics, Cambridge University Press, Cambridge, 1992.
\bibitem{zykov} A.A. Zykov, On some properties of linear complexes, Amer. Math. Soc. Transl. (1952) no. 79.
\end{thebibliography}
\end{document}